\documentclass{article}

\usepackage{amsmath} 
\usepackage{amssymb}

\title{\LARGE \bf Variable structure control for parabolic evolution 
                  equations
}

\author{Laura Levaggi%
\thanks{L. Levaggi is with Department of Mathematics, 
        University of Genova, Via Dodecaneso 35 - 16146 Genova, Italy
        {\tt\small levaggi@dima.unige.it}}%
}

\begin{document}


\newtheorem{lem}{Lemma}[section]
\newtheorem{prop}{Proposition}[section]
\newtheorem{teo}{Theorem}[section]
\newtheorem{cor}{Corollary}[section]
\newtheorem{hip}{Hypothesis}[section]


\newtheorem{defi}{Definition}[section]
\newtheorem{rem}{Remark}[section]
\newtheorem{ex}{Example}[section]


\def\beq{\begin{equation}}
\def\eeq{\end{equation}}
\newcommand {\edim} {\hfill\rule{0mm}{0mm}\hfill$\Box$\medskip\newline}
\newcommand {\eedim} {\hfill\rule{0mm}{0mm}\hfill$\Box$}
\newcommand {\bdim} {{\sc{Proof. }}}


\def\B{{\cal B}}
\def\L{{\cal L}}
\def\M{{\cal M}}
\def\dom{{\rm dom}\,}
\def\D{{\cal D}}
\def\R{{\cal R}}
\def\cB{\overline{\B}}
\def\inte{\bigcap_{\e>0}}
\def\cco{\overline{\rm co}\,}
\def\co{{\rm co}\,}
\def\sub{\subset}
\def\nsub{\sub\!\!\!\!\!/\;}
\def\span{{\rm span\ }}
\def\m{\!\raisebox{.2ex}{$\smallsetminus$}}
\def\sd{\oplus}
\def\ov{\overline}
\def\dim{{\rm dim}\,}
\def\und{\underline}


\def\l{\lambda}
\def\d{\delta}
\def\t{\tau}
\def\e{\varepsilon}
\def\w{\omega}
\def\r{\rho}
\def\s{\sigma}
\def\a{\alpha}
\def\b{\beta}
\def\f{\varphi}
\def\p{\psi}
\def\te{\theta}
\def\g{\gamma}
\def\O{\Omega}
\def\G{\Gamma}
\def\S{\Sigma}


\def\larr{\leftarrow}
\def\rarr{\rightarrow}
\def\weak{\rightharpoonup}
\def\dec{\raisebox{.2ex}{$~\scriptstyle {\!\searrow}~$}}
\def\ra{\rangle}
\def\la{\langle}
\def\lsim{\raisebox{-1ex}{$~\stackrel{\textstyle <}{\sim}~$}}
\def\gsim{\raisebox{-1ex}{$~\stackrel{\textstyle >}{\sim}~$}}


\newcommand{\er}{\relax\ifmmode I\!\!R\else$I\!\!R$\fi}
\newcommand{\en}{\relax\ifmmode I\!\!N\else$I\!\!N$\fi} 
\newcommand{\ec}{\relax\ifmmode I\!\!\!\!C\else$I\!\!\!\!C$\fi}
\newcommand{\qu}{\relax\ifmmode I\!\!\!\!Q\else$I\!\!\!\!Q$\fi}

\def\Re{\R\!{\rm e}\,}


\newcommand {\diffeq} [2] { \left\{ \begin{array} {l} #1 \\ #2 \end{array} \right. }
\newcommand {\pde} [6] { \begin{array}{ll} #1 & #2 \\ 
   #3 & #4 \\ #5 & #6 \end{array} }


\def\ueq{u_{{\rm eq}}}
\def\sign{{\rm sign}\,}
\def\invcb{(CB)^{-1}}
\def\ta{\tilde{A}}
\def\tK{\tilde{K}}
\newcommand {\parder} [3] {\frac{\partial^{#1} #2}{\partial #3^{#1}}}

\def\A{{\cal A}}
\def\C{{\cal C}}
\newcommand {\vct} [2] { \left( \begin{array} {c} #1 \\ #2 \end{array} \right) }

\maketitle
\date{}

\begin{abstract}

In this paper it is considered a class of infinite-dimensional control 
systems in a variational setting. By using a Faedo-Galerkin method, a 
sequence of approximating finite dimensional controlled differential 
equations is defined. On each of these systems a variable structure 
control is applied to constrain the motion on a specified surface. 
Under some growth assumptions the convergence of these approximations 
to an ideal sliding state for the infinite-dimensional system is shown.
Results are then applied to the Neumann boundary control of a parabolic 
evolution equation. 
\end{abstract}
%
\section{Introduction} \label{sec1}
Variable structure control methods and in particular sliding mode controls, 
are by now recognised as classical tools for the regulation of systems 
governed by ordinary differential equations in a finite dimensional 
setting. For an overview of the finite-dimensional theory see \cite{Utk92}.

While being easy to design, they possess attractive properties of 
robustness and insensitivity with respect to disturbances and unmodelled 
dynamics. These characteristics are all the more important when dealing 
with infinite-dimensional systems. In many control applications such as 
heat transfer processes, chemical processes, flexible manipulators the 
state evolution is governed by a partial differential equation. The 
complexity of these plants results in models having significant degrees 
of uncertainty. Thus motivated, recent research has been devoted to the 
extension of sliding mode control and therefore the use of discontinuous
feedback laws, to the infinite-dimensional setting. 
While earlier works \cite{OrlUtkARC82,OrlUtkAUTO87,OrlUtkAMCS98} were 
confined to some special classes of systems, at present both theory and 
application of sliding mode control have been extended to a rather 
general setting \cite{OrlTAC00,OrlDocTAC02,OrlLouChrIJC04,LevDIE02,
LevEJC02,LevCC04}. In particular in \cite{OrlTAC00} the key concept of 
equivalent control is introduced in a general Hilbert space framework 
for evolution equations governed by unbounded linear operators that 
generate $C_0$-semigroups. Also it is shown that, under some stability
assumptions, the ideal sliding can be uniformly approximated by ``real"
motions evolving in a boundary layer of the sliding manifold, thus 
ensuring the validity of the method for application purposes. The 
relationship between the equivalent control method and generalised 
solutions of infinite-dimensional systems with discontinuous 
right-hand side is presented in \cite{LevDIE02,LevEJC02}. 

All the results in the above cited literature only take into consideration 
distributed control systems, i.e. they deal with bounded input operators.
In this paper we make a first attempt to consider the extension of 
sliding modes to a class of boundary control problems in a general setting. 
To the author's knowledge there exist only a few results in this direction 
in the linear case \cite{DraUtkIJC92,DraOzgSV94}, where by application of 
integral transformations the problem is reduced to the control of a 
finite-dimensional differential-difference equation. Our approach goes 
instead in the direction of \cite{ZolB89}. In Section \ref{sec2} we define 
the general abstract variational framework in which we set up our control 
problem. In particular, the main assumptions we make on the operator governing 
the evolution, are weak continuity and coerciveness, so that both linear and 
non-linear operators are comprised in this setting. In Section \ref{sec3}
we present our main result: a Faedo-Galerkin method is used to construct 
a sequence of finite-dimensional approximations of the given problem. On
each of these the standard variable structure control theory of \cite{Utk92} 
can be applied. We then assume that for each approximation a control law is 
chosen to constrain the evolution in a boundary layer of a given sliding
manifold and study the limit as the dimensions diverge. We show that, under 
some growth assumption on the norm of these controls, a limit motion exists, 
which satisfies the sliding condition. Then, in Section \ref{sec4} we 
apply the obtained results to the Neumann boundary control of a heat equation.  
%
\section{Abstract setting and problem statement} \label{sec2}
In this paper we are going to consider a class of parabolic partial 
differential equations with controllers acting on the boundary. In 
particular we will study the case of Neumann boundary conditions and 
finite dimensional control space. Also, we suppose that a manifold 
$S$ is given, on which we want to restrict the motion of our system.  
We then analyse the problem of the existence of an admissible control 
law for which this ideal sliding motion is possible.  
\begin{ex}  \label{ex1}
Before going into the details of the precise abstract setting of the 
problem, we show an example of application to give an idea of the family 
of systems we intend to study. 

Let $\O$ be a bounded, open subset of $\er^n$ 
with smooth boundary $\G$, $T>0$ and $\Delta$ be the laplacian differential 
operator on $\er^n$. Consider the following evolution equation
\beq  \label{pdeinex}
  \pde {\parder{}{Q}{t}(t,x)=\Delta Q(t,x)+q(x)Q(t,x)} {t\in (0,T),\; x\in \O}
       {\parder{}{Q}{\nu}(t,\s)= u(t)g(\s)} {t\in (0,T),\; \s\in \G}
       {Q(0,x) = Q_0(x)} {x\in \O.}
\eeq
Here $Q:[0,T]\times \O\rarr \er$ represents the evolution of the ``state
vector", $u:[0,T]\rarr \er$ is a scalar control law, $g:\G\rarr \er$ 
and $q:\er^n\rarr \er$ is bounded. This equation represents a model of 
heat conduction with both diffusion and heat generation (if $q$ is 
nonnegative). 
Now for $\g:\O\rarr\er$ we can define (informally) a sliding surface $S$ 
as the set of functions $f:\O\rarr \er$ such that 
\[
  \int_\O f(x) \g(x)\, dx = 0
\]
In this case a sliding motion $Q(t,x)$ on $S$ would satisfy
\[
  \int_\O Q(t,x) \g(x)\, dx = 0, \quad t>0
\]
\end{ex}
\subsection{Variational formulation} \label{subsec2}
The setting of the abstract problem follows \cite{Lio69,Lio71,LioMag72}: 
let $V$ be a separable, reflexive Banach space, $H$ be a Hilbert space, 
$V\subset H$ with continuous injection. The space $H$ is identified with 
its dual, while we denote by $V'$ the dual space of $V$, so that we have 
\[
  V\subset H \subset V'.
\]
For $u_1$, $u_2\in H$ the scalar product in $H$ will be denoted by 
$(u_1,u_2)$ and the derived norm by $|u_i|$. We will denote by $\|\cdot\|$ 
the norm in $V$ and by $\|\cdot\|_*$ that in $V'$. The dual pairing 
between the two spaces will be written as $\la \cdot,\cdot \ra$.  
Also, we will assume that on $V$ it is defined a semi-norm $[\cdot]$ such 
that 
\beq  \label{vnorm}
  [v]+\l|v|\geq \b\|v\|, \quad \forall v\in V,\quad\mbox{for some }
  \l,\b >0.
\eeq
It is assumed that all the above (infinite-dimensional) spaces are 
real vector spaces; results can be extended to the complex case with 
the necessary modifications. For any $T>0$ we can define the following 
spaces of vector-valued functions: 
\[
  L^2(0,T;V)  =  \{f:[0,T]\rarr V: 
                \int_0^T\|f(t)\|^2 dt<+\infty\} 
\]
\[
  L^\infty(0,T;H) =  \{f:[0,T]\rarr H: 
               \sup_{t\in [0,T]}|f(t)| <+\infty\}.  
\]
The space $L^2(0,T;V')$ can be defined analogously. Also, it is possible
to define on these spaces a concept of derivative, in a distributional
sense (see i. e. \cite{Lio71} Chapter III). The following 
result \cite{LioMag72} will be useful in the sequel.
\begin{teo}\label{teocont}
  Let 
  \[
    W(0,T) = \left\{ f\in L^2(0,T;V)\,:\, \frac {df}{dt} \in L^2(0,T;V')
       \right \}.
  \]
  All functions in $W(0,T)$ are, after eventual modification on a null
  measure set, continuous from $[0,T]$ in $H$, i.e. $ W(0,T)\sub
  C^0(0,T;H)$.
\end{teo}
\vskip .3cm  \indent
For $t\in (0,T)$ let $A(t):V\rarr V'$ be an operator satisfying the 
following assumptions: 
\begin{itemize}
  \item for all $v,w\in V$ the map 
  \beq \label{meas}
     t\mapsto \la A(t)v,w\ra  \mbox{ is measurable;} 
  \eeq
  \item for all $t$ and any $ u,v,\w \in V$ the map
  \beq \label{hemi}
    \a\mapsto \la A(t)(u+\a v),w\ra  
    \mbox{ is continuous;} 
  \eeq
  \item there exist constants $c_1>0$, $c_2\geq 0$ such that 
  \beq \label{bound}
    \|A(t)v\|_*\leq c_1\|v\|+c_2, \quad \forall v\in V;
  \eeq  
  \item there exist constants $\a>0$ and $\nu\in \er$ such that
  \beq \label{coercive}
    \la A(t)v,v\ra \geq \a [v]^2 + \nu\, |v|^2  \forall v\in V.
  \eeq
  \item $A(\cdot)$ is $2$-weakly continuous, i.e.
  \begin{multline}  \label{2weak}
     v_k\rarr v\mbox{ weakly in } W(0,T) \Longrightarrow \\
     A(\cdot)v_k(\cdot)\rarr A(\cdot)v(\cdot) \mbox{ weakly in }
     L^2(0,T;V').
  \end{multline}
\end{itemize}
%
Let $U\sub \er^m$ be closed and convex and let $f:[0,T]\times U
\rarr V'$ satisfy the following condition: there exists a constant $C>0$
such that for any $u:[0,T]\rarr U$, $u\in L^2(0,T)$ 
\beq   \label{contf}
  \int_0^T\|f(t,u(t))\|^2_*\leq C \|u\|^2_2 ,
\eeq
where $\|u\|_2$ is the usual $L_2$-norm (it will always be understood
that control laws $u$ take values in $U$, so that we will write 
$u\in L^2(0,T)$ instead of $L^2(0,T;U)$). 

We are now ready to write the abstract evolution equation we are 
going to study. The evolution of the system will be given by a 
vector-valued function $y\in W(0,T)$ satisfying the following 
abstract Cauchy problem
\beq  \label{abspro}
  \diffeq { \frac {dy}{dt} + A(t)y(t) = f(t,u(t))\quad {\rm q.o.}\,t }
          {y(0)=y_0,}
\eeq
with $u\in L^2(0,T)$ and for some $y_0\in H$ (by Theorem \ref{teocont} 
this makes sense). The differential equation above as to be understood 
as an equality in the dual space $V'$, i.e. setting 
\beq   \label{defa}
  a(t;v,w) = \la A(t)v,w\ra,\quad t>0,\;v,w, \in V
\eeq
and in view of Theorem \ref{teocont}, the differential problem 
(\ref{abspro}) is equivalent to the following variational formulation
\beq \label{varpro}
  \diffeq { \frac{d}{dt}(y(t),v) + a(t;y(t),v) = \la f(t,u(t)),v\ra
           \;\; \forall v\in V, }
         {y(0)=y_0}
\eeq 
Existence and uniqueness results of the solution of such equations, 
under our assumptions, can be found in \cite{Lio69} under 
monotonicity assumptions and in \cite{Lio71,LioMag72} for the linear
case. 
%
\begin{ex}  \label{ex1-1}
  Let us see how Example \ref{ex1} fits into this framework. Let 
  $H=L^2(\O)$ and 
  \[
    V=H^1(\O)=\left \{f\in H:\,\parder{}{f}{x_i} \in H\,\,i=1,\ldots,n
    \right \}.
  \]  
  On $V$ we set $[v]^2=|\nabla v|^2$ and $\|v\|^2= [v]^2+|v|^2$.
  Let $v\in V$ be arbitrary; by scalar multiplication and using 
  Green's formula one finds that the solution $Q$ of (\ref{pdeinex}) 
  has to satisfy
  \begin{eqnarray*}
    \frac d{dt}(Q(t,\cdot),v) & = & \int_{\O} \Delta Q(t,x) 
        v(x)\,dx +(q\,Q(t,\cdot),v)\\
    & = & -\int_{\O} \nabla_x Q(t,v)\cdot \nabla v(x)\,dx \\
    &   & +\int_\G u(t)g(\s) \, v(\s)\,d\s +(q\,Q(t,\cdot),v).
  \end{eqnarray*}
  Therefore setting $y_0=Q_0$ and $y(t)=Q(t,\cdot)$ we get the 
  (autonomous) variational formulation of our abstract setting in 
  the form (\ref{varpro}) with
  \beq  \label{ainex}
    a(v,w) = (\nabla v,\nabla w)-(qv,w)
  \eeq
  and 
  \beq  \label{finex}
    \la f(t,u),v\ra = \int_\G ug(\s) \, v(\s)\,d\s.
  \eeq
  Now (\ref{hemi}) and (\ref{bound}) are easily verified and
  (\ref{coercive}) follows from
  \[
    a(v,v) = [v]^2 -(qv,v)\geq [v]^2-(\sup_{\O} q)\,|v|^2.
  \] 
  Also, the operator $A$ defined as $\la Av,w\ra:=a(v,w)$ is 
  linear and bounded, therefore it is $2$-weakly continuous and 
  we have (\ref{2weak}).

  Moreover, on $V$ the trace operator $\t$ of restriction of a 
  function to the boundary of $\O$ is well defined \cite{LioMag72}.
  The range of $\t$ is the Banach space $Z=H^{1/2}(\G)$ and $\g$ is 
  continuous from $V$ onto $H^{1/2}(\G)$. Therefore $f$ is well 
  defined for any $g$ in the dual of $H^{1/2}(\G)$, hence for 
  example for all $g\in L^2(\G)$ and obviously satisfies (\ref{contf})
  with $C=\|g\|_{L^2(\G)}\,\|\t\|_{\L(V,Z)}$. 
\end{ex}
%
\section{Main results}  \label{sec3}
In this section we introduce the concept of sliding surface for the 
control problem (\ref{varpro}) and show how sliding motions can be 
defined in this context.

Assume we are working in the framework set up in Section \ref{sec2}.
Thanks to separability, there exists a countable basis for $V$, so 
that it is possible to define a family $\{V_k\}_{k\in\en}$ of finite 
dimensional subspaces of $V$ 
\[
   V_k=\span \{v_{1,k},\ldots,v_{N_k,k}\}
\]
such that
\[
  V_k\sub V_{k+1},\quad \mathop{\bigcup_{k\in\en}} V_k=V.
\]
Then it is possible to define approximate solutions of (\ref{varpro})
by projecting on the subspaces $V_k$, using the standard Faedo-Galerkin 
method. We thus define the following family of variational problems:
find $y_k:[0,T]\rarr V_k$ such that
\beq   \label{varprok}
  \diffeq { \frac{d}{dt}(y_k(t),v) + a(t;y_k(t),v) = \la f(t,u_k(t)),v\ra
           \; \forall v\in V_k, }
         {y_k(0)=y_{0,k}}
\eeq
with $y_{0,k}\in V_k$ for all $k$ and a sequence $\{u_k\}$ in $L^2(0,T)$.
Note that, since $V_k$ has dimension $N_k$, the above problem can be 
written as an ordinary differential equation. In fact, since 
$y_k(t)\in V_k$, there exists a vector $\xi_k(t) \in \er^{N_k}$ 
such that 
\[
  y_k(t)=\sum_{i=1}^{N_K} (\xi_k(t))_i\,v_{i,k}.
\]
The differential equation in (\ref{varprok}) is satisfied for all 
$v\in V_k$ iff it is valid for every element of the basis of $V_k$. 
Therefore, if
\[
  y_{0,k}=\sum_{i=1}^{N_K} (\xi_{0,k})_i \,v_{i,k} ,
\]
\[
  f_k(t)=(\,\la f(t,u_k(t)),v_{i,k}\ra\,)_{i=1, \ldots,N_k}, 
\] 
and 
\begin{eqnarray*}
  A_k(t)=(a^{(k)}_{ij}(t))_{i,j=1,\ldots,N_k}, & \; &
   a^{k}_{ij}(t)=a(t;v_{i,k},v_{j,k}),\\
  M_k=(m^{(k)}_{ij}(t))_{i,j=1,\ldots,N_k}, & \; &
   m^{k}_{ij}(t)=(v_{i,k},v_{j,k}),
\end{eqnarray*}
the differential problem (\ref{varprok}) is equivalent to the 
following ordinary Cauchy problem
\beq  \label{odek}
 \diffeq{M_k\,\dot{\xi}_k(t)+A_k\,\xi_k(t)=f_k(t)}
        {\xi_k(0)=\xi_{0,k}.}
\eeq
We now prove a convergence result for the approximations $y_k$, under
some conditions on the controls sequence $\{u_k\}$. 
\begin{teo}  \label{teoconv}
  Let the assumptions in Section \ref{sec2} be satisfied and $\{u_k\}$ 
  be a sequence in $L^2(0,T)$. Let $y_k$ be the solution of 
  (\ref{varprok}) and suppose that $y_{0,k}\rarr y_0$ in $H$ for 
  $k\rarr +\infty$. Suppose moreover that the following condition on the 
  growth of the control norms is satisfied
  \beq   \label{congrow}
    \|u_k\|_{L^2(0,t)}^2\leq M\int_0^t\,|y_k(s)|^2\,ds+N,
     \; t\leq T
  \eeq
  for some non-negative constants $M$ and $N$ and that $f$ is 
  the following weak continuity assumption 
  \begin{multline}  \label{weakcontf}
    u_k\rarr u^* \mbox{ weakly in } L^2(0,T) \mbox{ then } \\
    f(\cdot,u_k(\cdot))\rarr f(\cdot,u(\cdot))\mbox{ weakly in }
    L^2(0,T;V')
  \end{multline} 
  Then there exist a
  control law $u^*\in L^2(0,T)$ and a function $y^*\in W(0,T)$ 
  verifying (\ref{varpro}), such that, for some subsequence, 
  \begin{eqnarray*}
    & & y_k\rarr y^* \mbox{ weakly in } W(0,T) \\
    & & y_k\rarr y^* \mbox{ weakly* in } L^\infty(0,T;H)\\
    & & u_k\rarr u^* \mbox{ weakly in } L^2(0,T).
  \end{eqnarray*}
\end{teo}
\vskip .3cm
\bdim
  Writing (\ref{varprok}) for $v=y_k(t)$ we get
  \[
    (\dot{y}_k(t),y_k(t))+a(t;y_k(t),y_k(t))=\la f(t,u_k(t)),y_k(t)\ra.
  \]
  As the first term on the left is in fact the time derivative of 
  $|y_k(t)|^2/2$, integrating the above identity we have
  \begin{multline*}
    \frac 12 |y_k(t)|^2 + \int_0^t a(t;y_k(s),y_k(s))\,ds = \\
       \frac 12 |y_k(0)|^2 +\int_0^t \la f(t,u_k(s)),y_k(s)\ra\,ds.
  \end{multline*}
  By (\ref{coercive}), (\ref{contf}) and (\ref{vnorm}) we obtain the 
  following inequality
  \begin{multline*}
    \frac 12  |y_k(t)|^2 + \a\int_0^t [y_k(s)]^2 ds \leq \\ 
       \frac 12 |y_k(0)|^2 - \nu \int_0^t |y_k(s)|^2 ds \\
        +c\,\|u_k\|_2 \,\left( \int_0^t [y_k(s)]^2 \,ds
            +\int_0^t [y_k(s)]^2\,ds\right)^{1/2}
  \end{multline*}
  for some constant $c>0$. Consider now for $x\geq 0$ the function 
  $h(x)=(\a x)/2-c\sqrt{x}$. It is easy to show that it has minimum 
  for $x=(c/\a)^2$, therefore $c\sqrt{x}\leq (\a x+c^2/\a)/2$, thus
  \begin{multline*} 
    \frac 12  |y_k(t)|^2 + \frac \a 2\int_0^t [y_k(s)]^2 ds \leq \\
    \frac 12 |y_k(0)|^2 + \left(\frac\a 2+|v|\right) 
    \int_0^t |y_k(s)|^2 ds +\frac {c^2}{2\a}\|u_k\|_2^2.
  \end{multline*}
  Now, since by hypothesis $|y_{0,k}-y_0|$ tends to zero, the 
  term $|y_k(0)|^2$ is bounded. Moreover by (\ref{congrow}) 
  \beq \label{ineqen}
    |y_k(t)|^2 +\a \int_0^t [y_k(s)]^2 ds \leq 
    c_1 + c_2 \int_0^t |y_k(s)|^2 \,ds
  \eeq
  for some constants $c_1,c_2>0$. Since $\a>0$ we get
  \[
    |y_k(t)|^2\leq c_1+c_2\int_0^t |y_k(s)|^2\,ds
  \]
  Therefore, by Gronwall's lemma we obtain for some constant 
  $K>0$
  \beq \label{boundH}
    \|y_k\|_{L^{\infty}(0,T;H)} = \sup_{t\in [0,T]}
                      |y_k(t)|\leq K
  \eeq
  therefore from (\ref{ineqen}) we also have
  \[
    \int_0^T[y_k(s)]^2\,ds \leq {\rm const}
  \]
  and lastly, using (\ref{vnorm}) and (\ref{bound})
  \begin{eqnarray*}
   \|y_k\|_{L^2(0,T;V)} =\left (\int_0^T\|y_k(s)\|^2\,ds \right ) 
   \leq {\rm const}, \\
   \int_0^T\|A(t)y_k(t)\|_*\,dt\leq {\rm const}. 
  \end{eqnarray*} 
  Since spheres are weakly compact in both $L^2(0,T;V)$ and 
  $L^2(0,T;V')$, weakly* compact in $L^{\infty}(0,T;H)$, we can 
  extract a 
  subsequence of $\{y_k\}$ (which for simplicity we still denote 
  by $\{y_k\}$) converging to some $y^*\in L^2(0,T;V)\cap 
  L^{\infty}(0,T;H)$ for both the weak topology of $L^2(0,T;V)$ 
  and the weak* topology of $L^{\infty}(0,T;H)$ and such that
  $Ay_k$ weakly converges to some $\eta$ in $L^2(0,T;V')$. By 
  (\ref{congrow}) we also have that $\|u_k\|_2$ is bounded, thus 
  eventually passing to a further subsequence, there exists 
  $u^*\in L^2(0,T)$ such that $u_k$ converges to $u^*$ weakly in 
  $L^2(0,T)$. Also, by (\ref{weakcontf}) we can proceed as in 
  the proof of Theorem 1.1, p. 159 of \cite{Lio69} to conclude 
  that 
  \[
    \diffeq {\frac d{dt} y^*(t)+\eta(t) = f(t,u^*(t))}
            {y(0)=y_0.}
  \]
  Also, by a standard argument (see i.e. \cite{ZolB89}, Theorem 3)
  one can prove that $\dot{y}_k\rarr \dot{y}^*$ weakly in 
  $L^2(0,T;V')$, i.e. $y_k\rarr y^*$ weakly in $W(0,T)$. Thus, 
  by (\ref{2weak}) $\eta(t)=A(t)y^*(t)$ and the proof is complete.
\edim
Having achieved the above convergence result, we introduce as in 
\cite{ZolB89}, a set $D$ which can be either $V$ or a sufficiently 
large open subset of $H$ and a mapping $s:D\rarr \er^m$ continuously 
Fr\'echet differentiable on $D$. The sliding 
surface $S$ we consider is defined as $S=\{y\in D:\,s(y)=0\}$. 
Proceeding as in \cite{ZolB89}, by slightly modifying proofs, it is 
possible to prove the following
\begin{cor}  \label{corsli}
  Let the assumptions of Theorem \ref{teoconv} hold. Let $z_k(t)
  =s(y_k(t))$ and assume that one of the following is satisfied:
  \begin{itemize}
    \item [(1)] $D=V$, $s$ is affine and $z_k\rarr 0$ uniformly
                in $t$;
    \item [(2)] $\B_H(0,K)\sub D$, $V$ is compactly embedded 
                in $H$ (here $\B_H$ denotes a ball in $H$, while 
                $K$ is defined in (\ref{boundH}) above) and 
                $z_k(t)\rarr 0$ for almost every $t\in [0,T]$.
  \end{itemize}
  Then the limit motion $y^*$ of Theorem \ref{teoconv} belongs to 
  the sliding manifold $S$.
\end{cor}
\begin{rem}
  Note that by (\ref{odek}) every $y_k$ solves a finite-dimensional
  problem, thus for the approximate solutions all results of the 
  classical theory of variable structure systems and sliding mode 
  control of \cite{Utk92} are valid. Therefore existence results
  for system motions satisfying the requirements in Corollary 
  \ref{corsli} and design methods to achieve them are available. 
  See also the discussion of existence under relaxed hypotheses 
  developed in \cite{ZolB89}.
\end{rem}
\addtolength{\textheight}{-2.8cm}
%
\section{An application}  \label{sec4}
In this Section we show an application of the obtained results on 
the control problem introduced in Example \ref{ex1}. We have already 
proved (see Example \ref{ex1-1}) that this partial differential 
equation with Neumann control fits in the abstract setting of Section
\ref{sec2}. It is also easy to prove that for $f$ as in (\ref{finex})
the condition (\ref{weakcontf}) is satisfied. In fact, if $u_k\rarr u$
weakly in $L^2(0,T)$, for any $\f\in L^2(0,T;V)$ we have
\begin{multline*}
  \int_0^T \la\, [f(t,u_k(t))- f(t,u(t))]\,,\,\f(t)\,\ra \,dt = \\ 
     \int_0^T[u_k(t)-u(t)]\, \int_\G g(\s) \, \f(t)(\s)\,d\s\;dt
\end{multline*}
which converges to zero since by H\"older's inequality and 
continuity of the trace operator on $V$
\begin{multline*}
  \int_0^T \left( \int_\G |g(\s)| \, |\f(t)(\s)|\,d\s 
  \right)^2dt \leq \\
  \|g\|_{L^2(\G)}^2 \int_0^T\|\f(t)\|^2\,dt < +\infty.
\end{multline*} 

We then set $s:H\rarr \er$, $s(x)=(x,\g)$ and $S=\ker S$. 
For convenience we suppose that the chosen bases of the subspaces 
$V_k$ are orthonormal, so that the matrix $M_k$ in (\ref{odek}) is 
the identity (this is not restrictive since in the general case 
$M_k$ is symmetric, positive definite and a linear change of 
coordinates is sufficient to reconduct this problem to the 
orthonormal one). Then, setting $g_k=(\, (g,\t v_{i,k})
_{L^2(\G)}\,)_{i=1, \ldots,N_k}$, (\ref{odek}) can be rewritten 
as
\[
  \diffeq{\dot{\xi}_k(t)+A_k\,\xi_k(t)=u_k(t)g_k}
        {\xi_k(0)=\xi_{0,k}.}
\]
Then $z_k(t)=s(y_k(t))=(y_k(t),\g)=\g_k^T\xi_k(t)$, with 
$\g_k=(\, (v_{i,k},\g)\,)_{i=1, \ldots,N_k}$. Let $V(t)=
z_k^2(t)/2$; then 
\[
  \dot{V}(t) = z_k(t)\,\dot{z}_k(t) = z_k(t)\,[\g_k^T\,
           (-A_k\xi_k(t) +u_k(t) g_k)\,].
\]
By standard finite dimensional theory \cite{Utk92} 
a sliding mode exists on $S_k=\{x\in\er^{N_k}:\,\g_k^Tx=0\}$ if
$\g_k^Tg_k \neq 0$. Also, in this case, setting 
\[
  u_k(t) = -U(t)\; \frac{\sign(z_k(t))}{\g_k^Tg_k}
\]
with $U(t)> |\g_k^TA_k\xi_k(t)|$ the sliding 
surface is globally attractive and reached in finite time. 
Moreover, if $\d_k>0$ and $|s(y_k(0))|<\d_k$ the control 
\[
  u_k(t) = -\frac{U(t)}{\g_k^Tg_k}\; \frac{z_k(t)}{|z_k(t)|+\d_k}
\]
constrains the motion of the system in a $\d_k$-boundary layer 
of $S_k$. Let us now consider the term $\g_k^TA_k\xi_k(t)$; 
since we assumed that the basis of $V_k$ is orthonormal, we have
\[
  \g_k^TA_k\xi_k(t) = a(y_k(t),P_k\,\g),
\]
where $P_k:V\rarr V_k$ is the projection on $V_k$. Likewise we 
have
\[
  \g_k^Tg_k = \int_\G g(\s)\,(P_k\,\g)(\s)\,ds.
\]
Thus, if for example $\g\in V$ and 
\[
  \int_\G g(\s)\,\g(\s)\,d\s\neq 0,
\]
since $P_k\,\g\rarr \g$ in $V$, there exists $K$ such that 
$\g_k^Tg_k\neq 0$ for all $k\geq K$. In order to apply Theorem
\ref{teoconv} we also have to show that (\ref{congrow}) holds. 
Recalling that 
\[
  a(y_k(t),P_k\,\g) =(\nabla y_k(t),\nabla P_k\,\g)
  -(qy_k(t),P_k\,\g)
\]
we just have to show that, at least for suitable $\g$-s, the first 
term can be estimated using $|y_k(t)|$. Proceeding formally, by 
Green's formula we have
\begin{multline*}
  (\nabla y_k(t),\nabla P_k\,\g)= -(y_k(t),\Delta(P_k\,\g)) \\
   + \int_\G y_k(t)(\s)\,\parder{}{}{\nu}(P_k\,\g) (\s)\,d\s.
\end{multline*}
Thus (\ref{congrow}) can be satisfied if sufficiently regular 
decompositions $\{V_k\}$ of $H^1(\O)$ are chosen and if the 
function $\g$ satisfies $\parder{}{}{\nu}P_k\,\g=0$, at least 
on some subsequence. For example this is true if $\g\in V_N$ 
for some $N$ and $\parder{}{}{\nu}\g=0$.  
\begin{rem}
  In this paper we have chosen a variational setting for our 
  problem, by which we can encompass also some non-linear partial 
  differential equations. For the linear case, another common
  abstract setting involves semigroup theory. In the above 
  example our operator $A:V\rarr V'$ could be, in some sense,  
  substituted by $\A:\D(\A)\sub H\rarr H$, with 
  \[
    \D(\A) = \left \{ y\in H^2(\O):\, \parder{}{y}{\nu} =0 
        \right \},\\
    \quad \A y = \Delta y+qy.
  \]
  Note that the last condition on $\g$ above is related to 
  ``$\g\in\D(\A^*)$", which is frequently encountered in the 
  literature on output control of infinite-dimensional systems.
\end{rem}
\begin{rem}
  In many applications the function $z(t)$ of the example represents
  the system's output. The modulus of the control law we have chosen 
  depends on the whole state norm, which could be unavailable for 
  measurement. In \cite{OrlLouChrIJC04} observers are designed to 
  overcome this difficulty in the case of distributed control. It 
  would be interesting to study their application to this case also.
\end{rem} 
%
\section{Conclusions and future work}
In this paper we have analysed the convergence behaviour of finite 
dimensional Faedo-Galerkin approximations of a class of variational 
problems, when sliding motions are taken into consideration. We have
thus shown that, under some growth hypothesis on the norms of the 
controls, a sliding motion exists.

This is a first attempt to extend variable structure control 
to boundary control problems for infinite-dimensional systems 
and much work has still to be done in this area. Apart from the need 
to extend these results to different boundary control problems, it 
would be interesting to study how these results are related to a
notion of equivalent control, which has already by introduced in 
the infinite-dimensional setting and to approximability of ideal 
sliding motions by real ones. 

\bibliographystyle{plain}
\bibliography{bdcon}

\end{document}